\documentclass[a4paper,UKenglish,cleveref, autoref, thm-restate,pdfa]{lipics-v2021}

\pdfoutput=1 
\hideLIPIcs  

\title{Formalization of $p$-adic $L$-functions in Lean 3} 

\titlerunning{$p$-adic $L$-functions} 

\author{Ashvni Narayanan}{London School of Geometry and Number Theory, Imperial College London}{a.narayanan20@imperial.ac.uk}{https://orcid.org/0000-0003-2777-4228}{EPSRC Grant EP/S021590/1 (UK)}

\authorrunning{A. Narayanan} 

\Copyright{Ashvni Narayanan} 



\supplement{Copies of the source files relevant to this paper are available in a separate repository.}
\supplementdetails[swhid={}]{Software}{https://github.com/laughinggas/p-adic-L-functions} 

\acknowledgements{The author is supported by EPSRC. The author would like to thank her PhD supervisor Prof Kevin Buzzard for several helpful insights. 
She would also like to thank Dr Filippo A. E. Nuccio for helpful conversations that helped give shape to the project. She is very grateful 
to the entire Lean community for their timely help and consistent support.}

\nolinenumbers 

\usepackage[utf8]{inputenc}
\usepackage{amsmath}
\usepackage{amssymb}
\usepackage{amsfonts}
\usepackage{hyperref}
\usepackage{makecell}
\usepackage{xcolor}
\usepackage{xspace}
\usepackage{listings}

\definecolor{keywordcolor}{rgb}{0.7, 0.1, 0.1}   
\definecolor{commentcolor}{rgb}{0.4, 0.4, 0.4}   
\definecolor{symbolcolor}{rgb}{0.4, 0.4, 0.4}    
\definecolor{sortcolor}{rgb}{0.1, 0.5, 0.1}      

\newcommand{\lean}[1]{\texttt{#1}\xspace} 

\newcommand*{\Cl}{\mathcal{C}\kern-.075em l}

\newcommand{\N}{\mathbb{N}}

\newcommand{\Q}{\mathbb{Q}}
\newcommand{\Z}{\mathbb{Z}}

\DeclareUnicodeCharacter{2208}{\ensuremath{\in}}
\DeclareUnicodeCharacter{21D1}{\ensuremath{\Uparrow}}
\DeclareUnicodeCharacter{2194}{\ensuremath{\longleftrightarrow}}
\DeclareUnicodeCharacter{1D4DD}{\ensuremath{\mathcal{N}}}
\DeclareUnicodeCharacter{22C3}{\ensuremath{\bigcup}}
\DeclareUnicodeCharacter{03BB}{\ensuremath{\lambda}}
\DeclareUnicodeCharacter{03B1}{\ensuremath{\alpha}}
\DeclareUnicodeCharacter{03C3}{\ensuremath{\sigma}}
\DeclareUnicodeCharacter{03C7}{\ensuremath{\chi}}
\DeclareUnicodeCharacter{2081}{\ensuremath{_1}}
\DeclareUnicodeCharacter{2082}{\ensuremath{_2}}
\DeclareUnicodeCharacter{2090}{\ensuremath{_a}}
\DeclareUnicodeCharacter{2097}{\ensuremath{_l}}
\DeclareUnicodeCharacter{211A}{\ensuremath{\Q}}
\DeclareUnicodeCharacter{2124}{\ensuremath{\Z}}
\DeclareUnicodeCharacter{2115}{\ensuremath{\N}}
\DeclareUnicodeCharacter{2211}{\ensuremath{\sum}}
\DeclareUnicodeCharacter{2264}{\ensuremath{\l}}
\DeclareUnicodeCharacter{2286}{\ensuremath{\subseteq}}
\DeclareUnicodeCharacter{22A4}{\ensuremath{\top}}
\DeclareUnicodeCharacter{22A5}{\ensuremath{\bot}}
\DeclareUnicodeCharacter{2243}{\ensuremath{\simeq}}
\DeclareUnicodeCharacter{2080}{\ensuremath{_0}}
\DeclareUnicodeCharacter{03C8}{\ensuremath{\psi}}
\DeclareUnicodeCharacter{2225}{\ensuremath{\lVert}}
\DeclareUnicodeCharacter{2A06}{\ensuremath{\bigsqcup}}
\DeclareUnicodeCharacter{02E3}{\ensuremath{^{\times}}}
\DeclareUnicodeCharacter{03B5}{\ensuremath{\varepsilon}}
\DeclareUnicodeCharacter{2203}{\ensuremath{\exists}}
\DeclareUnicodeCharacter{2200}{\ensuremath{\forall}}
\DeclareUnicodeCharacter{2265}{\ensuremath{\ge}}
\DeclareUnicodeCharacter{1DA0}{\ensuremath{^{\texttt{f}}}}

\begin{document}

\maketitle

\begin{abstract}
The Euler--Riemann zeta function is a largely studied numbertheoretic object, and the birthplace of several conjectures, 
such as the Riemann Hypothesis. Different approaches are used to study it, including $p$-adic analysis : 
deriving information from $p$-adic zeta functions. A generalized version of $p$-adic zeta functions (Riemann zeta function) 
are $p$-adic $L$-functions (resp. Dirichlet $L$-functions). This paper describes formalization of $p$-adic $L$-functions 
in an interactive theorem prover Lean 3. Kubota--Leopoldt $p$-adic $L$-functions are meromorphic functions emerging from the 
special values they take at negative integers in terms of generalized Bernoulli numbers. They also take twisted values of 
the Dirichlet $L$-function at negative integers. This work has never been done before in any theorem prover. 
Our work is done with the support of \lean{mathlib} 3, one of Lean's mathematical libraries. It required formalization of a 
lot of associated topics, such as Dirichlet characters, Bernoulli polynomials etc. We formalize these first, 
then the definition of a $p$-adic $L$-function in terms of an integral with respect to the Bernoulli measure, 
proving that they take the required values at negative integers.
\end{abstract}

\section{Introduction}
We are working on formalizing mathematics in an interactive theorem prover called Lean. 
Formal verification involves the use of logical and computational methods to establish 
claims that are expressed in precise mathematical terms \cite{TPIL}. Lean is a powerful tool 
that facilitates formalization of a system of mathematics supported by a basic set of axioms. There is a large mathematical library of theorems verified by Lean called \lean{mathlib}, maintained by a community of computer scientists and mathematicians. One can then 
formally verify proofs of new theorems dependent on preexisting theorems in \lean{mathlib}. 
\lean{mathlib} contained 100579 theorems(as of early October 2022). It would be impossible to construct 
such a vast library without a highly collaborative spirit and a communal decentralized effort, one of 
Lean's best features. 

$p$-adic $L$-functions are a well studied numbertheoretic object. They were initially constructed 
by Kubota and Leopoldt in \cite{KL}. Their motivation was to construct a meromorphic function that helps study 
the Kummer congruence for Bernoulli numbers, and gives information regarding $p$-adic class numbers. 
As a result, these functions take twisted values of the Dirichlet $L$-function at negative integers, and 
are also related to the generalized Bernoulli numbers and the $p$-adic zeta function. There are several different ways of 
constructing $p$-adic $L$-functions, we refer to the constructions
given in Chapter 12 of \cite{cyc}. As a result, one needs to build a lot of background (in the maximum possible generality) 
before embarking on the main goal.

It is difficult to explain all the mathematical terms used here, we attempt to describe as many as possible. To that effect, 
a basic knowledge of algebra is assumed. Since \lean{mathlib} works in utmost generality, one often finds that the terminology 
used is less common. Thanks to the community's endeavour to maintain adequate documentation, we have added links which serve as 
explanations wherever possible. When clear, we will explicitly skip writing hypotheses in the code, since these can get quite long.

We give a mathematical overview in this section, then discuss background in Section \ref{section2}, define Dirichlet characters in Section \ref{dirchar}, 
introduce generalized Bernoulli numbers in Section \ref{ber}, construct the $p$-adic $L$-function in Section \ref{section3}, and evaluate it at negative integers in Section \ref{section4}, finishing with a summary in Section \ref{section5}. 
\subsection{Mathematical overview}
We give a brief overview of the mathematics formalized in this project. 
$L$-functions are a fundamental object, appearing almost everywhere in modern 
number theory. The Dirichlet $L$-function associated to a Dirichlet character $\chi$ is given by  
$$ L(s, \chi) = \sum_{n = 1}^{\infty} \frac{\chi (n)}{n^s} = \prod_{p \text{ prime}} \frac{1}{1 - \chi (p) p^{-s}}$$
where $s$ is a complex variable with $Re(s) > 1$. This can be analytically extended 
to the entire complex plane, with a simple pole at $s = 1$ when $\chi = 1$. Note also 
that $L(s, 1)$ is the same as the Riemann zeta function. Moreover, it is known that 
$L(1 - n, \chi) = - \frac{B_{n, \chi}}{n}$, where $B_{n, \chi}$ are the generalized 
Bernoulli numbers. 

In this paper, we construct, for an integer prime $p$, a $p$-adic analogue of $L(s, \chi)$, 
called the Kubota--Leopoldt $p$-adic $L$-function, denoted $L_p(s, \chi)$. This is generally done by continuously extending the function 
$L_p(1 - n, \chi) := (1 - \chi (p) p^{n - 1}) L(1 - n, \chi)$ to the complete $p$-adic space 
$\mathbb{C}_p$. In fact, $L_p(s, 1)$ is analytic except for a pole at $s = 1$ with residue 
$1 - \frac{1}{p}$ (Theorem 5.11, \cite{cyc}). 

Formalization of the $p$-adic $L$-functions via analytic continuation was hard, since $\mathbb{C}_p$ did not exist 
in \lean{mathlib} at the time. Following \cite{cyc}, 
we instead define it in terms of an ``integral'' with respect to the Bernoulli measure. We explain these terms below. 

A profinite space is a compact, Hausdorff and totally disconnected space. The $p$-adic integers 
$\mathbb{Z}_p$, which are the completion of the integers $\mathbb{Z}$ with respect to the 
valuation $\nu_p(p^{\alpha} \prod_{p_i \ne p} p_i^{\alpha_i}) = \alpha$ are a profinite space. 
One may also think of them as the inverse limit of the discrete topological spaces 
$\mathbb{Z} / p^{n} \mathbb{Z}$, that is, $\mathbb{Z}_p = \projlim_{n} \mathbb{Z} / p^{n} \mathbb{Z}$. 

Locally constant functions are those for which the preimage of any set is open. Given a profinite space $X$ and a normed ring $R$, 
one can show that the locally constant functions from $X$ to $R$ (denoted $LC(X, R)$) are dense in the space of continuous 
functions from $X$ to $R$ (denoted $C(X, R)$).

Given an abelian group $A$, a distribution is defined to be an $A$-linear map from $LC(X, A)$ 
to $A$. A measure $\phi$ is defined to be a bounded distribution, that is, $\forall f \in LC(X, R)$, 
$\exists K > 0$ such that $|| \phi (f) || \le K ||f|| $, where $||f|| = \sup_{x \in X} || f(x)||$. 
An example of a measure is the Bernoulli measure. Given a natural number $d$ coprime to $p$ and a clopen set $U_{n, a}$ of 
$\mathbb{Z}/ d \mathbb{Z} \times \mathbb{Z}_p$, the characteristic function $\chi_{n, a}$ 
(defined to be 1 on $U_{n, a}$ and 0 otherwise) is a locally constant function. 
Given a natural number $c$ that is coprime to $d$ and $p$, we then define the Bernoulli measure $E_c$ by :
$$ E_c(\chi_{n, a}) := \bigg \{ \frac{a}{d p^{n + 1}} \bigg \} - c \bigg \{ \frac{c^{-1} a}{d p^{n + 1}} \bigg \} + \frac{c - 1}{2} $$

Given a measure $\mu$, the integral with respect to $\mu$ is 
$\int f d\mu := \mu (f)$ for any locally constant function $f$, and 
extending this definition to $C(X, R)$. In fact, this is an $R$-linear map. 

Finally, the $p$-adic $L$-function is defined to be an integral with respect to the Bernoulli 
measure. The characterizing property of the $p$-adic $L$-function is its evaluation at negative integers : 
$$ L_p (1 - n, \chi) = -(1 - \chi \omega^{-n}(p)p^{n - 1}) \frac{B_{n, \chi \omega^{-n}}}{n} $$
for $n \ge 1$. When defined as an integral, additional work is needed to prove this. 

Our contributions to this theory include a formalized definition of the $p$-adic $L$-function in generality, 
taking values in a normed complete non-Archimedean $\mathbb{Q}_p$-algebra, instead of just $\mathbb{C}_p$. Further, it takes as 
input continuous monoid homomorphisms, also known as elements of the weight space. We have also developed an extensive theory for 
Dirichlet characters, Bernoulli numbers and polynomials, generalized Bernoulli numbers, properties of $p$-adic integers and modular arithmetic, 
making substantial contributions to the \lean{number\_theory} section of \lean{mathlib}. 
We use non-traditional methods to define and prove classical results, often choosing to work with those which are easier to formalize, 
later proving their equivalence to the original. 

\subsection{Lean and \lean{mathlib}}
Lean 3 is a functional programming language and interactive theorem prover based on
dependent type theory. This project is based on Lean’s 
mathematical library \lean{mathlib 3}, which is characterized by
its decentralized nature with over 300 contributors. Thus, it is impossible to cite every author who contributed a piece of code that we used.

We assume the reader is familiar with structures such as \lean{def}, \lean{abbreviation}, \lean{lemma}, \lean{theorem}, which are used constantly. 
An important property of Lean is its typeclass inference system - 
Lean ``remembers'' properties given to a \lean{structure} or \lean{class} embedded in an \lean{instance} structure. This is explained in detail in \cite{mathlib}. 
We shall also use several tactics in proofs, such as \lean{rw}, \lean{apply}, \lean{conv} and \lean{refine}
\footnote{\url{https://leanprover-community.github.io/mathlib_docs/tactics.html} has a full list of tactics in Lean}.
\section{Preliminaries}
\label{section2}
\subsection{Filters and convergence}
None of our mathematical proofs require filters on paper, however, we find that working with them makes 
formalizing our proofs significantly less cumbersome. Due to the efforts of 
Johannes Hölzl, Jeremy Avigad, Patrick Massot and several others, we have a vast API for filters in Lean. 
We shall not delve into the details of what a filter is, but instead explain how they are used to formalize convergence and limits. 
\newline For a sequence of functions $(f_n)_{n \in \mathbb{N}}$, the expression $\lim_{n \to \infty} f_n(x) = a$ is 
represented as :
\begin{lstlisting}
tendsto (λ n : ℕ, f_n) filter.at_top ( 𝓝 a)
\end{lstlisting}
Here, \href{https://leanprover-community.github.io/mathlib_docs/order/filter/at_top_bot.html#filter.at_top}{\lean{filter.at\_top}} 
(for the naturals) is a filter on $\mathbb{N}$ generated by the collection of sets $\{ b \mid a \leq b \}$ 
for all $a \in \mathbb{N}$. 
The following lemma is particulary useful :
\begin{lstlisting}
/-- If f₁ and f₂ are equal almost everywhere, then f₁ converges if and only if f₂ converges. -/
lemma filter.tendsto_congr' {α : Type} {β : Type u_1} {f₁ f₂ : α → β} 
  {l₁ : filter α} {l₂ : filter β} (h : f₁ =ᶠ[l₁] f₂) : 
  tendsto f₁ l₁ l₂ ↔ tendsto f₂ l₁ l₂
\end{lstlisting}
This lemma shows that sequences that are the same after finitely many 
elements have the same limit. 
Given two sequences \lean{f₁} and \lean{f₂} (thought of as functions from $\mathbb{N}$), \newline \lean{f₁ =ᶠ[at\_top] f₂ ↔ ∃ (a : ℕ), ∀ (b : ℕ), b ≥ a, f₁ b = f₂ b}. 

An equivalent condition to convergence on metric spaces is : 
\begin{lstlisting}
lemma metric.tendsto_at_top : ∀ {α : Type u_1} {β : Type} 
  [pseudo_metric_space α] [nonempty β] [semilattice_sup β] 
  {u : β → α} {a : α} : 
  tendsto u at_top ( 𝓝 a) ↔ ∀ (ε : ℝ) (h : ε > 0), 
  (∃ (N : β), ∀ (n : β), n ≥ N → ‖ u n - a ‖ < ε)
\end{lstlisting}
Thus, in order to prove lemmas about convergence, one can either choose to continue doing computations in the \lean{tendsto} framework, 
or prove normed inequalities. Working with the former really simplified 
calculations. 
As an example, suppose we want to prove the convergence of the sequence $g$ given by $g(0) = g(2) = 1$ and $g(n) = 3 f(n)$, 
where $f$ is a convergent sequence. This is a one-line proof using \href{https://leanprover-community.github.io/mathlib_docs/order/filter/basic.html#filter.tendsto_congr%27}{\lean{filter.tendsto\_congr'}}. 
Using the above lemma, one must obtain $N$ corresponding to $\varepsilon / 3$, and also prove that $0 < \varepsilon / 3$. 
With more complex expressions, this gets computationally difficult to handle. 

Hence, we try to avoid using \href{https://leanprover-community.github.io/mathlib_docs/topology/metric_space/basic.html#metric.tendsto_at_top}{\lean{metric.tendsto\_at\_top}} 
when possible. The only cases where it is used is when direct inequalities 
need to be dealt with; this happens precisely when the non-Archimedean condition on $R$ is used. Thus, this is a good 
indicator of where the non-Archimedean condition is needed.
\subsection{Modular arithmetic and units}
Some fundamental objects with which we shall work throughout are the finite spaces $\mathbb{Z}/n \mathbb{Z}$. Note that proving properties for \href{https://leanprover-community.github.io/mathlib_docs/data/zmod/defs.html#zmod}{\lean{zmod n}} 
is equivalent to proving them for any finite cyclic group. Given a positive $n \in \mathbb{N}$, 
$\mathbb{Z}/n \mathbb{Z}$ is the same as \href{https://leanprover-community.github.io/mathlib_docs/init/data/fin/basic.html#fin}{\lean{fin n}}(and $\mathbb{Z}$ for $n = 0$), 
the set of natural numbers upto $n$. It is also the set of equivalence classes obtained via the relation 
on $\mathbb{Z}$ : $a \sim b \iff n \mid a - b$. It has a natural group structure, and is given the discrete topology, 
making it a topological group. Some maps used constantly include \newline \href{https://leanprover-community.github.io/mathlib_docs/data/zmod/basic.html#zmod.val}{\lean{val:zmod n → ℕ}}, 
which takes any element to its smallest nonnegative reprentative less than \lean{n}; 
and \href{https://leanprover-community.github.io/mathlib_docs/data/zmod/basic.html#zmod.cast_hom}{\lean{cast\_hom:zmod n → R}}, a coercion to a ring, 
obtained by composing the canonical coercion with \lean{val}. If \lean{R} has characteristic 
dividing \lean{n}, the map is a ring homomorphism. Given coprime naturals $m, n$, an important equivalence is 
\newline \href{https://leanprover-community.github.io/mathlib_docs/data/zmod/basic.html#zmod.chinese_remainder}{\lean{chinese\_remainder:zmod (m * n) ≃+* zmod m × zmod n}}. 
About 45 additional lemmas were required, which have been put in a separate file, 
\href[]{https://github.com/laughinggas/p-adic-L-functions/blob/main/src/zmod/properties.lean}{\lean{zmod/properties.lean}}. 

Every \href{https://leanprover-community.github.io/mathlib_docs/algebra/group/defs.html#monoid}{\lean{monoid}} \lean{M} has an associated space of 
invertible elements or units, denoted \newline \href{https://leanprover-community.github.io/mathlib_docs/algebra/group/units.html#units}{\lean{units M}} or \lean{Mˣ}. 
We use the map \lean{units.coe\_hom:Mˣ → M} to identify a unit in its parent space frequently. 
Given a \href{https://leanprover-community.github.io/mathlib_docs/algebra/hom/group.html#monoid_hom}{\lean{monoid\_hom}} (abbreviated as \lean{→*}) \lean{R →* S} for monoids \lean{R} and \lean{S}, 
one can obtain a homomorphism \lean{Rˣ →* Sˣ} by \href{https://leanprover-community.github.io/mathlib_docs/algebra/hom/units.html#units.map}{\lean{units.map}}. 

\section{Dirichlet characters and the Teichmüller character}
\label{dirchar}
An important task was to formalize Dirichlet characters, an integral part of the definition of the $p$-adic $L$-function. 
Dirichlet characters are often not found to be defined in this technical manner. Another addition is the definition of 
Dirichlet characters of level and conductor 0. The words character and Dirichlet character are used interchangeably. 

Dirichlet characters are usually defined as group homomorphisms from $\mathbb{Z}/n \mathbb{Z}^{\times}$ to $\mathbb{C}^{\times}$ for some natural number $n$. 
A lot of properties traditionally known for groups hold more generally and are defined in greater generality in \lean{mathlib} for monoids. In the same spirit, 
we define Dirichlet characters to be monoid homomorphisms on any \href{https://leanprover-community.github.io/mathlib_docs/algebra/group/defs.html#monoid}{\lean{monoid}} :
\begin{lstlisting}
abbreviation dirichlet_character (R : Type*) [monoid R] (n : ℕ) := 
  (zmod n)ˣ →* R ˣ 
/-- The level of a Dirichlet character. -/
abbreviation lev {R : Type*} [monoid R] {n : ℕ} 
  (χ : dirichlet_character R n) : ℕ := n
\end{lstlisting}
If we gave the definition of Dirichlet characters a \lean{def} structure, 
\lean{dirichlet\_character} would become a \lean{Type} distinct from \lean{(zmod n)ˣ →* Rˣ}, 
making compositions with \lean{monoid\_hom} complicated; hence we used \lean{abbreviation} instead. 
Note that the linter returns an extra unused argument warning (for \lean{χ}) for the latter definition. 

Given a Dirichlet character $\chi$, \lean{asso\_dirichlet\_character χ} returns a monoid homomorphism from $\mathbb{Z}/n \mathbb{Z}$ 
to $R$, which is $\chi$ on the units and 0 otherwise. 
\begin{lstlisting}
noncomputable abbreviation asso_dirichlet_character {R : Type*} 
  [monoid_with_zero R] {n : ℕ} (χ : dirichlet_character R n) : 
  zmod n →* R := { to_fun := 
  function.extend (units.coe_hom (zmod n)) ((units.coe_hom R) ∘ χ) 0, ..}
\end{lstlisting}
Lean requires us to tag this definition \lean{noncomputable}, since we are producing data from an existential statement, 
\lean{classical.some} (appearing in \lean{function.extend}), which has no computational content (see Chapter 11 of \cite{TPIL}).
One would like to shift between compatible Dirichlet characters of different levels. For this, we construct the following tools : 
\begin{lstlisting}
/-- Extends the Dirichlet character χ of level n to level m, for n | m. -/
def change_level {m : ℕ} (hm : n | m) : 
  dirichlet_character R n →* dirichlet_character R m :=
{ to_fun := λ ψ, ψ.comp (units.map (zmod.cast_hom hm (zmod n))), .. }
\end{lstlisting}
\newpage
\begin{lstlisting}
/-- χ₀ of level d factors through χ of level n if d | n and 
    χ₀ = χ ∘ (zmod n → zmod d). -/
structure factors_through (d : ℕ) : Prop :=
(dvd : d | n) 
(ind_char : ∃ χ₀ : dirichlet_character R d, χ = χ₀.change_level dvd)
\end{lstlisting}
The notions of primitivity and conductor of a Dirichlet character follow easily : 
\begin{lstlisting}
/-- The set of numbers for which a Dirichlet character is periodic. -/
def conductor_set : set ℕ := {x : ℕ | χ.factors_through x}
/-- The minimum natural number n for which a character is periodic. -/
noncomputable def conductor : ℕ := Inf (conductor_set χ)
/-- A character is primitive if its level is equal to its conductor. -/
def is_primitive : Prop := χ.conductor = n
/-- The primitive character associated to a Dirichlet character. -/
noncomputable def asso_primitive_character : dirichlet_character R χ.conductor := classical.some (χ.factors_through_conductor).ind_char
\end{lstlisting}
Here, \href{https://leanprover-community.github.io/mathlib_docs/init/classical.html#classical.some}{\lean{classical.some}} 
makes an arbitrary choice of an element from a nonempty space, and \href{https://leanprover-community.github.io/mathlib_docs/init/classical.html#classical.some_spec}{\lean{classical.some\_spec}} 
lists down properties of this element coming from the space. 

When $a = b$, while \lean{dirichlet\_character R a} and \lean{dirichlet\_character R b} are ``mathematically'' equal, Lean does 
not think of them as the same type. This gets complicated when additional layers, such as \lean{change\_level} are added to the equation. 
A general method to resolve such problems is by using the tactic \href{https://leanprover-community.github.io/mathlib_docs/tactics.html#subst}{\lean{subst}}, 
which would substitute $a$ with $b$; however, that failed. Instead, we used the concept of heterogeneous equality (\href{https://leanprover-community.github.io/mathlib_docs/init/core.html#heq}{\lean{heq}}, or \lean{==}) 
to deal with this. The tactic \href{https://leanprover-community.github.io/mathlib_docs/tactics.html#congr'}{\lean{congr'}} helped reduce to expressions of heterogeneous equality, 
which were then solved with the help of lemmas such as :
\begin{lstlisting}
lemma change_level_heq {a b : ℕ} {S : Type*} [comm_monoid_with_zero S] 
  (χ : dirichlet_character S a) (h : a = b) : 
  change_level (show a | b, from by {rw h}) χ == χ
\end{lstlisting}
This states that, for $a = b$, changing the level of a Dirichlet character of level $a$ to $b$ is heterogeneously equal to itself. \newline
Traditionally only for primitive characters, our definition of multiplication of characters 
extends to any two characters. This takes as input characters \lean{χ₁} and \lean{χ₂} of levels \lean{n} and \lean{m} respectively, and returns the primitive 
character associated to $\chi_1' \chi_2'$, where $\chi_1'$ and $\chi_2'$ are obtained by changing the levels 
of \lean{χ₁} and \lean{χ₂} to \lean{lcm n m}. 
\begin{lstlisting}
noncomputable def mul {m n : ℕ} (χ₁ : dirichlet_character R n) 
  (χ₂ : dirichlet_character R m) :=
asso_primitive_character(change_level χ₁ (dvd_lcm_left n m) * 
  change_level χ₂ (dvd_lcm_right n m))
\end{lstlisting} 
This multiplication is not trivially commutative or associative, with respect to this definition. 

We need the notion of odd and even characters. A character $\chi$ is odd if $\chi (-1) = -1$, and 
even if $\chi (-1) = 1$. For a commutative ring, any character is either odd or even : 
\begin{lstlisting}
lemma is_odd_or_is_even {S : Type*} [comm_ring S] [no_zero_divisors S] 
  {m : ℕ} (ψ : dirichlet_character S m) : ψ.is_odd ∨ ψ.is_even
\end{lstlisting}
\subsection{Teichmüller character}
The initial effort was to formalize the definition of the Teichmüller character (denoted $\omega$) directly. 
However, it was discovered that Witt vectors, and in particular Teichmüller lifts had previously 
been added to \lean{mathlib} by Johan Commelin and Robert Lewis. This reiterates the importance 
of the collaborative spirit of Lean, and of making definitions in the correct generality. 

It is beyond the scope of this text to define Witt vectors and do it justice. We refer interested readers to 
Section 2.4 of \cite{witt}. For a commutative ring $R$ and a prime number $p$, one can obtain a ring of Witt vectors $\mathbb{W}(R)$.
When we take $R = \mathbb{Z}/p \mathbb{Z}$, we get 
\begin{lstlisting}
def equiv : $\mathbb{W}$ (zmod p) ≃+* ℤ_[p]
\end{lstlisting}
One also obtains the Teichmüller lift $R \to \mathbb{W} (R)$. Given $r \in R$, 
the 0-th coefficient is $r$, and the other coefficients are 0. This map is a 
multiplicative monoid homomorphism and is denoted 
\href{https://leanprover-community.github.io/mathlib_docs/ring_theory/witt_vector/teichmuller.html#witt_vector.teichmuller}{\lean{teichmuller}}. 

Combining this with the previous two definitions, we obtain our definition of the Teichmüller character : 
\begin{lstlisting}
noncomputable abbreviation teichmuller_character_mod_p (p : ℕ) 
  [fact (nat.prime p)] : dirichlet_character ℤ_[p] p := units.map 
  (((witt_vector.equiv p).to_monoid_hom).comp (witt_vector.teichmuller p))
\end{lstlisting}
We use \lean{[fact p.prime]} to make the primality of $p$ an instance. This map takes \newline $x \in \mathbb{Z}/p\mathbb{Z}^{\times}$ to a root of unity $y \in \mathbb{Z}_p$ such that $y \equiv x (\texttt{mod p})$. 
Often we view this as taking values in a $\mathbb{Q}_p$-algebra $R$, by composing it 
with \lean{algebra\_map ℚ\_[p] R}, which identifies elements of $\mathbb{Q}_p$ 
in $R$. Since we mostly deal with $\omega ^{-1}$ taking values on $R^{\times}$, we define this as \lean{teichmuller\_character\_mod\_p'}. \newline
We proved properties of Teichmüller characters in \href{https://github.com/laughinggas/p-adic-L-functions/blob/main/src/dirichlet_character/teichmuller_character.lean}{\lean{teichmuller\_character.lean}}, 
such as, for odd primes $p$, the Teichmüller character is odd, and \lean{1} otherwise : 
\begin{lstlisting}
lemma eval_neg_one (hp : 2 < p) : teichmuller_character_mod_p p (-1) = -1
\end{lstlisting}

\section{Bernoulli polynomials and the generalized Bernoulli number}
\label{ber}
The Bernoulli numbers $B_n'$ are generating functions given by $\sum B_n'\frac{t^n}{n!}=\frac{t}{e^{t} - 1}$. They appear 
in the computation of sums of powers of naturals, $\sum_n n^k$. Note that several authors think of Bernoulli numbers $B_n$ 
to be defined as $\sum B_n\frac{t^n}{n!}=\frac{t}{1-e^{-t}}$. The difference between these two is : $B_n = (-1)^n B_n'$, 
with $B_1' = - \frac{1}{2}$.
A reformulation gives :
$$ B_n' = 1 - \sum_{k = 0}^{n - 1} {n \choose k} \frac{B_k'}{n - k + 1} $$
In \lean{mathlib}, $B_n'$ was already defined (by Johan Commelin) as above. However, we needed $B_n$, which we then defined as :
\begin{lstlisting}
def bernoulli (n : ℕ) : ℚ := (-1)^n * bernoulli' n
\end{lstlisting}
The Bernoulli polynomials, denoted $B_n(X)$, a generalization of the Bernoulli numbers,
are generating functions $ \sum_{n = 0}^{\infty} B_n(X) \frac{t^n}{n!} = \frac{t e^{tX}}{e^t - 1} $. 
This gives :
$$ B_n (X) = \sum_{i = 0}^n {n \choose i} B_i X^{n - i} $$
We defined the Bernoulli polynomials as : 
\begin{lstlisting}
def polynomial.bernoulli (n : ℕ) : polynomial ℚ :=
  ∑ i in range (n + 1), monomial (n - i) ((bernoulli i) * (choose n i))
\end{lstlisting}
Here, \href{https://leanprover-community.github.io/mathlib_docs/data/polynomial/basic.html#polynomial.monomial}{\lean{monomial n a}} 
translates to $a X^n$, and \lean{∑ i in s, f i} translates to $\sum_{i \in s} f(i)$, for a \lean{finset} (or finite set) \lean{s}. A small aspect of this naming convention is that if the namespaces for Bernoulli numbers and polynomials are both open 
(which is often the case), in order to use the Bernoulli numbers, one needs to use \lean{\_root\_.bernoulli}. 
We shall use them interchangeably here, when the context is clear. 



An important fact is, $\forall n$, $(n + 1)X^{n} = \sum_{k = 0}^{n} {n \choose k} B_k(X)$ :
\begin{lstlisting}
theorem sum_bernoulli (n : ℕ) : monomial n (n + 1 : ℚ) = 
  ∑ k in range (n + 1), ((n + 1).choose k : ℚ) • bernoulli k 
\end{lstlisting}
These proofs are relatively straightforward. Most of this work is now part of \lean{mathlib}, 
and has been used to give a formalized proof of \href{https://leanprover-community.github.io/mathlib_docs/number_theory/bernoulli.html\#sum_range_pow}{Faulhaber's theorem}.
\subsection{Generalized Bernoulli numbers}
Generalized Bernoulli numbers are integral to our work, since these are related to the special values of $p$-adic $L$-functions and Dirichlet $L$-functions. 
Given a primitive Dirichlet character $\chi$ of conductor $f$, the generalized Bernoulli numbers are defined as (section 4.1, \cite{cyc}) 
$ \sum_{n = 0}^{\infty} B_{n,\chi} \frac{t^n}{n!} = \sum_{a = 1}^f \frac{\chi(a)t e^{at}}{e^{ft} - 1} $. 
For any multiple $F$ of $f$, Proposition 4.1 of \cite{cyc} gives us : 
$$ B_{n, \chi} = F^{n - 1} \sum_{a = 1}^{F} \chi (a) B_n \bigg( \frac{a}{F} \bigg) $$
This is much easier to work with, so we use this as our definition instead, taking $F = f$ :
\begin{lstlisting}
def general_bernoulli_number {S : Type*} [comm_semiring S] [algebra ℚ S] 
  {n : ℕ} (ψ : dirichlet_character S n) (m : ℕ) : S :=
  (algebra_map ℚ S ((ψ.conductor)^(m - 1 : ℤ))) * 
  ∑ a in finset.range ψ.conductor,
  asso_dirichlet_character (asso_primitive_character ψ) a.succ * 
  algebra_map ℚ S ((bernoulli m).eval (a.succ / ψ.conductor : ℚ))
\end{lstlisting}
Contrary to the traditional definition, this is for all characters, and $\psi$ takes values 
in any commutative $\mathbb{Q}$-algebra, instead of $\mathbb{C}$. One had to also explicitly mention that \lean{m - 1} must be taken to have type \lean{ℤ}, since Lean would otherwise infer 
it to have type \lean{ℕ}, which might have caused errors (subtraction on $\mathbb{N}$ 
and $\mathbb{Z}$ are different).

\subsection{A special property of generalized Bernoulli numbers}
An important property of these numbers is : 
\begin{theorem}\label{thm1}
Let $\chi$ be an even Dirichlet character of level $dp^m$ for $d$ coprime to the odd prime $p$, with $m$ positive. 
Suppose $R$ is a nontrivial commutative non-Archimedean normed $\mathbb{Q}_p$-algebra with no zero divisors. 
For $k > 1$, 
$$ \lim_{n \to \infty} \frac{1}{dp^{n}} \sum_{0 < a < dp^{n} ; (a, dp) = 1} \chi \omega^{-k} (a) a^{k} = 
(1 - \chi \omega^{-k} (p) p^{k-1}) B_{k, \chi \omega^{-k}} $$
\end{theorem} 
Instead of giving \lean{R} a non-Archimedean structure(which did not exist in \lean{mathlib} when this project began), 
we give as input its consequences, conditions \lean{na} and \lean{na'}. This is formulated in Lean as :
\begin{lstlisting}
theorem lim_even_character' (na' : ∀ (n : ℕ) (f : (zmod n)ˣ → R), 
  ∥∑ i : (zmod n)ˣ, f i ∥ ≤ ⨆ (i : (zmod n)ˣ), ∥f i ∥)
(na : ∀ (n : ℕ) (f : ℕ → R), 
  ∥ ∑ (i : ℕ) in finset.range n, f i ∥ ≤ ⨆ (i : zmod n), ∥f i.val ∥) : 
  tendsto (λ (n : ℕ), (1 / ↑(d * p ^ n)) • 
  ∑ (i : ℕ) in finset.range (d * p ^ n),
  asso_dirichlet_character (χ.mul (teichmuller_character_mod_p' p R ^ k)) ↑i * ↑i ^ k) at_top ( 𝓝 (general_bernoulli_number 
  (χ.mul (teichmuller_character_mod_p' p R ^ k)) k))
\end{lstlisting}
The proof of this theorem follows from the proof in Lemma 7.11 of \cite{cyc}, a 
point of difference being that our theorem holds more generally for $R$ being a non-Archimedean normed commutative $\mathbb{Q}_p$-algebra with no zero divisors, 
instead of $\mathbb{C}_p$. Majorly, it equates the two sides modulo $p^n$ for a 
sufficiently large $n$, and uses the fact that 
\begin{theorem}\label{thm2}
  $$ \lim_{n \to \infty} \frac{1}{dp^{n}} \sum_{0 < a < dp^{n} ; (a, dp) = 1} \chi \omega^{-m} (a) a^{m} = 0 $$
\end{theorem}
The formalization is very calculation intensive, and is a good example of a small proof on paper being magnified in Lean, 
because there are multiple coercions and arithmetic calculations to be dealt with. 
Unfortunately, tactics such as \href{https://leanprover-community.github.io/mathlib_docs/tactics.html#ring}{\lean{ring}} and 
\href{https://leanprover-community.github.io/mathlib_docs/tactics.html#simp}{\lean{simp}} that usually help with these fail here. 
It is translated in Lean as :
\begin{lstlisting}
lemma sum_even_character_tendsto_zero_of_units :
  tendsto (λ n, ∑ (i : (zmod (d * p^n))ˣ), ((asso_dirichlet_character
  (χ.mul (teichmuller_character_mod_p' p R^k))) i * i^(k - 1))) 
  at_top ( 𝓝 0) 
\end{lstlisting}
The proof of this theorem is in \href{https://github.com/laughinggas/p-adic-L-functions/blob/main/src/tendsto_zero_of_sum_even_char.lean}{\lean{tendsto\_zero\_of\_sum\_even\_char.lean}}.
\section{Construction of the $p$-adic $L$-function}
\label{section3}
\subsection{Density of locally constant functions}
For any compact Hausdorff totally disconnected space X and a commutative
normed ring A, we have proved that $LC(X, A)$ is a dense subset of $C(X, A)$. 
Formalizing this took about 500 lines of code (now in \lean{mathlib}), and is based on the fact that locally compact Hausdorff
totally disconnected spaces have a clopen basis :
\begin{lstlisting}
lemma loc_compact_Haus_tot_disc_of_zero_dim {H : Type*} [t2_space H] 
  [locally_compact_space H] [totally_disconnected_space H] : 
  is_topological_basis {s : set H | is_clopen s}
\end{lstlisting}
This turned out to be hard to formalize. Given a set $s$ of $H$, Lean gives a subset $V$ of $s$ 
the type \lean{V:set s}; however, Lean does not recognize $V$ as a subset of $H$. As a result, 
to use \lean{compact\_space s ↔ is\_compact (s:set H)}, one must construct \lean{V':set H} 
to be the image of $V$ under the closed embedding \lean{coe:s → H}. This process must be repeated each time a subset of $H$, 
which is also a topological subspace, is considered. Finally, it must be shown that all these coercions match up in the big 
topological space $H$.

\subsection{Clopen sets of the $p$-adic integers}
$\mathbb{Z}_p$ is a profinite space (as shown in section 2.4 of \cite{witt}). It is the inverse limit of finite
discrete topological spaces $\mathbb{Z}/p^n \mathbb{Z}$ for all $n$, and has a clopen basis of the
form $U_{a,n} := proj_n ^{-1} (a)$ for $a \in \mathbb{Z}/p^n \mathbb{Z}$, where $proj_n$ is the
canonical projection ring homomorphism \lean{to\_zmod\_pow n:$\mathbb{Z}$\_[p] →+* zmod (p $\hat{}$ n)}. 
We first define the collection of sets $(U_{a,n})_{a,n}$ :
\begin{lstlisting}
def clopen_basis : set (set ℤ_[p]) := 
  {x : set ℤ_[p] | ∃ (n : ℕ) (a : zmod (p^n)), 
    x = set.preimage (padic_int.to_zmod_pow n) {a} }
\end{lstlisting}
We show that \lean{clopen\_basis} forms a topological basis and that every element is
clopen :
\begin{lstlisting}
theorem clopen_basis_clopen : (clopen_basis p).is_topological_basis ∧ 
  ∀ x ∈ (clopen_basis p), is_clopen x
\end{lstlisting}
The mathematical proof is to show that for any $\epsilon$-ball, one can find $U_{a,n}$ inside it.
This is true because, given $n \in \mathbb{N}$ and $x \in \mathbb{Z} / p^n \mathbb{Z}$, the preimage of $x$ 
under \lean{to\_zmod\_pow n} is the same as the ball centered at $x$ (now considered as an element of $\mathbb{Z}_p$) with radius $p^{1 - n}$. 
The following lemmas prove useful :
\begin{lstlisting}
lemma appr_spec (n : ℕ) (x : ℤ_[p]) : 
  x - appr x n ∈ (ideal.span {p^n} : ideal ℤ_[p])
lemma has_coe_t_eq_coe (x : ℤ_[p]) (n : ℕ) : 
  (((appr x n) : zmod (p^n)) : ℤ_[p]) = ((appr x n) : ℤ_[p])
\end{lstlisting}
For \lean{x:ℤ\_[p]}, \lean{appr x n} is the smallest natural number in \lean{x (mod p\textasciicircum n)}. 
In the latter lemma, the RHS is a coercion of \lean{appr x n}, which has type \lean{$\mathbb{N}$}, to \lean{$\mathbb{Z}$\_p}. 
The LHS is a coercion of \lean{appr x n} to \lean{zmod (p $\hat{}$ n)} to \lean{$\mathbb{Z}$\_p}. 
This statement is not true in general, that is, given any natural number $n$, it is not true that the lift of $n$ to $\mathbb{Z}_p$ 
is the same as the composition of its lift to $\mathbb{Z}/p^n \mathbb{Z}$ and $\mathbb{Z}_p$. 
It works here because the coercion from $\mathbb{Z}/p^n \mathbb{Z}$ to $\mathbb{Z}_p$ is not the canonical lift.
It is a composition of a coercion from $\mathbb{Z}/p^n \mathbb{Z}$ to $\mathbb{N}$, which takes
$a \in \mathbb{Z}/p^n \mathbb{Z}$ to the smallest natural number in its
$\mathbb{Z}/p^n \mathbb{Z}$ equivalence class. 

One can similarly show that the sets $U_{b, a, n} := proj_1^{-1} (b) \times proj_{2,n} ^{-1} (a)$ form a clopen basis for 
$\mathbb{Z} / d \mathbb{Z} \times \mathbb{Z}_p$, where $proj_1$ is the first canonical projection on $b \in \mathbb{Z} / d \mathbb{Z}$ 
and $proj_{2,n}$ the composition of the second projection on $a \in \mathbb{Z}_p$ with $proj_n$ described above. We call this set 
\lean{clopen\_basis' p d}. Its properties are formalized in \href{https://github.com/laughinggas/p-adic-L-functions/blob/main/src/padic_int/clopen_properties.lean}{\lean{padic\_int.clopen\_properties.lean}}.

\subsection{Distributions and measures}
In this section, $X = \varprojlim_{i \in \mathbb{N}} X_i$ denotes a profinite space with $X_i$ finite and
projection maps $\pi_i : X \xrightarrow[]{} X_i$ and surjective maps
$\pi_{ij} : X_i \xrightarrow[]{} X_j$ for all $i \ge j$. Henceforth, we use $G$ to denote an abelian group,
$A$ for a commutative normed ring, $R$ for a commutative complete normed ring which is also a $\mathbb{Q}_p$-algebra, 
and $LC(X,Y)$ for the space of locally constant functions from $X$ to $Y$. 
We fix a prime $p$ and an integer $d$ such that $gcd(d, p) =1$. \newline
The topology on $C(X, A)$ comes from its normed group structure induced by the norm on $A$ :
$|| f - g || = sup_{x \in X} || f(x) - g(x) ||$. In fact, this topology is the same as the 
topology defined on bounded functions on $X$, since $X$ is a compact space. Since the API for bounded 
continuous functions on compact spaces was developed at around the same time (created by Oliver Nash), 
we used the existing lemmas such as 
\href{https://github.com/leanprover-community/mathlib/blob/32253a1a1071173b33dc7d6a218cf722c6feb514/src/topology/continuous_function/compact.lean#L47}{\lean{equiv\_bounded\_of\_compact}}.

A distribution (from Section 12.1 of \cite{cyc}) is a $G$-linear function $\phi : LC(X, G) \xrightarrow[]{} G$. 
This is already a \lean{Type}, hence we do not redefine it. 
Measures (not to be confused with measure theory measures) are bounded distributions :
\begin{lstlisting}
def measures := {φ : (locally_constant X A) →ₗ[A] A // ∃ K : ℝ, 0 < K ∧ 
    ∀ f : (locally_constant X A), $\lVert$φ f$\rVert$ ≤ K * $\lVert$inclusion X A f$\rVert$ }
\end{lstlisting}
The map \lean{inclusion} identifies the locally constant function \lean{f} as a continuous function. 
The boundedness of the distribution makes the measure continuous.
\subsection{The Bernoulli measure}
The Bernoulli measure is an essential measure. We make a choice of an integer $c$ with
$gcd(c,dp) = 1$, and $c^{-1}$ is an integer such that $c c^{-1} \equiv 1 \text{ mod } dp^{2n+1}$.
For a clopen set $U_{a,n}$, we define
$$ E_c (\chi_{U_{a,n}}) = E_{c,n} (a) = \bigg\{ \frac{a}{dp^{n + 1}} \bigg\} - c \bigg\{ \frac{c^{-1}a}{dp^{n + 1}} \bigg\} + \frac{c - 1}{2} $$
In Lean, this translates to (note that \lean{fract x} represents the fractional part of $x$) :
\begin{lstlisting}
def bernoulli_distribution := λ (n : ℕ) (a : (zmod (d * (p^n)))), 
  fract ((a : ℤ) / (d*p^(n + 1))) 
  - c * fract ((a : ℤ) / (c * (d*p^(n + 1)))) + (c - 1)/2
\end{lstlisting}
The original plan was to define a set of the form : 
\begin{lstlisting}
def bernoulli_measure (hc : c.gcd p = 1) :=
 {x : locally_constant (zmod d × ℤ_[p]) R →ₗ[R] R | ∀ (n : ℕ) 
  (a : zmod (d * (p^n))), x (char_fn R (clopen_from.is_clopen p d n a)) = 
  (algebra_map ℚ R) (E_c p d hc n a) }
\end{lstlisting}
and to show that it is nonempty. \lean{char\_fn} is a locally constant characteristic function on a clopen set ($1$ on the set and $0$ otherwise), 
taking as input the range of the function and the fact that the set is clopen. However, information is lost this way, since one then has to use 
\lean{classical.some} to extract the underlying measure. We use an elegant approach :
\begin{lstlisting}
/-- A sequence has the `is_eventually_constant` predicate if all the elements of the sequence are eventually the same. -/
def is_eventually_constant {α : Type*} (a : ℕ → α) : Prop := 
  { n | ∀ m, n ≤ m → a (nat.succ m) = a m }.nonempty
structure eventually_constant_seq {α : Type*} :=
(to_seq : ℕ → α) 
(is_eventually_const : is_eventually_constant to_seq)
\end{lstlisting}
Given a locally constant function $f$ from $\mathbb{Z}/d \mathbb{Z} \times \mathbb{Z}_p$ to $R$, 
we define the eventually constant sequence \lean{from\_loc\_const} : 
\begin{lstlisting}
noncomputable abbreviation from_loc_const : @eventually_constant_seq R :=
{ to_seq := λ (n : ℕ), 
    ∑ a in (zmod' (d * p^n) _), 
    f(a) • ((algebra_map ℚ_[p] R) (bernoulli_distribution p d c n a)),
  is_eventually_constant := _, }
\end{lstlisting} 
for all natural numbers $n$. \lean{zmod'} is the universal \lean{finset} of \lean{zmod}. 
We shall look into the proof of this sequence being eventually constant later. 

Given a locally constant function \lean{f:locally\_constant ((zmod d)ˣ × ℤ\_[p]ˣ) R}, 
an element of the set \lean{bernoulli\_measure} is given by : 
\begin{lstlisting}
  sequence_limit (from_loc_const p d R (loc_const_ind_fn _ p d f))
\end{lstlisting}
where \lean{loc\_const\_ind\_fn} is a locally constant function on $\mathbb{Z}/d \mathbb{Z} \times \mathbb{Z}_p$ 
that takes value $f$ on the units of the domain, and 0 otherwise. The linearity properties follow easily. Notice that \lean{bernoulli\_distribution} takes locally constant functions on $\mathbb{Z}/d \mathbb{Z} \times \mathbb{Z}_p$, 
while \lean{bernoulli\_measure} takes locally constant functions on $\mathbb{Z}/d \mathbb{Z}^* \times \mathbb{Z}_p^*$. This had to be done since our clopen basis was defined on 
$\mathbb{Z}/d \mathbb{Z} \times \mathbb{Z}_p$, and while it is easy to show the same results for the units on paper, it requires a bit of work in Lean.

We now prove that \lean{bernoulli\_measure} is indeed a measure, that is, it is bounded. The bound we choose is 
$K := 1 + \parallel c \parallel + \parallel \frac{c - 1}{2} \parallel$. The proof is as follows : let $\phi$ denote \lean{loc\_const\_ind\_fn}. 
We want $ \parallel E_c (\phi (f)) \parallel \le K \parallel f \parallel $. It suffices to prove this for $\chi_{n, a}$, because one can find an $n$ such that 
$\phi (f) = \sum_{a \in \mathbb{Z}/d \mathbb{Z} \times \mathbb{Z} /p^n \mathbb{Z}} \phi(f) (a) \dot{} \chi_{n,a}$ :
\begin{lstlisting}
lemma loc_const_eq_sum_char_fn (f : locally_constant ((zmod d) × ℤ_[p]) R) (hd : d.gcd p = 1) : ∃ n : ℕ, f = ∑ a in (finset.range (d * p^n)), 
      f(a) • char_fn R (clopen_from.is_clopen p d n a)
\end{lstlisting}
This proof is akin to proving that \lean{from\_loc\_const} is eventually constant, using discrete quotients. 
The discrete quotient on a topological space is given by an equivalence relation such 
that all equivalence classes are clopen : 
\begin{lstlisting}
structure (X : Type*) [topological_space X] discrete_quotient :=
(rel : X → X → Prop) 
(equiv : equivalence rel) 
(clopen : ∀ x, is_clopen (set_of (rel x)))
\end{lstlisting}
The last statement translates to, $\forall x \in X, \{ y | y \sim x \}$ is clopen. 
Given two discrete quotients $A$ and $B$, $A \le B$ means $\forall x,y \in X$, 
$x \sim_{A} y \implies x \sim_{B} y$. Any locally constant function induces a 
discrete quotient via its clopen fibers : 
\begin{lstlisting}
def locally_constant.discrete_quotient : discrete_quotient X := 
{ rel := λ a b, f b = f a, .. }
\end{lstlisting}
We now define a function : 
\begin{lstlisting}
/-- A discrete quotient induced by `to_zmod_pow`. -/
def discrete_quotient_of_to_zmod_pow : 
  ℕ → discrete_quotient (zmod d × ℤ_[p]) := 
  λ n, ⟨λ a b, to_zmod_pow n a.2 = to_zmod_pow n b.2 ∧ a.1 = b.1, _, _⟩
\end{lstlisting}
For $a = (a_1, a_2)$ and $b = (b_1, b_2)$ in $\mathbb{Z}/d \mathbb{Z} \times \mathbb{Z}_p$, this represents the relation 
$ a \sim b \iff a_2 (\text{mod } p^n) = b_2 (\text{mod } p^n) \wedge a_1 = b_1 $. 
Then, given a locally constant function $f$ on $\mathbb{Z}/d \mathbb{Z} \times \mathbb{Z}_p$, for $N$ large enough, the fibers of $f$ mod $p^N$ are contained in the basic clopen sets of $p^N$ :
\begin{lstlisting}
lemma le : ∃ N : ℕ, 
  discrete_quotient_of_to_zmod_pow p d N ≤ discrete_quotient f
\end{lstlisting} 
The proofs now follow from this fact : $\exists N, \forall m \ge N$,
$$ \sum_{a \in \mathbb{Z}/dp^{m + 1} \mathbb{Z}} f(a) E_{c,m + 1}(a) = \sum_{a \in \mathbb{Z}/dp^{m} \mathbb{Z}} f(a) E_{c,m}(a) $$
The required $N$ is \lean{classical.some (discrete\_quotient\_of\_to\_zmod\_pow.le f) + 1}. We also define the following : 
\begin{lstlisting}
/-- Set of all `b ∈ zmod (d * p^m)` such that `b = a mod (d * p^n)` for 
  `a ∈ zmod (d * p^n)`. -/
def equi_class (n m : ℕ) (a : zmod (d * p^n)) := 
  {b : zmod (d * p^m) | (b : zmod (d * p^n)) = a}
\end{lstlisting}
Then, we have the following lemma :
\begin{lstlisting}
lemma zmod'_succ_eq_bUnion : 
  zmod' (d*p^(m + 1)) = (zmod' (d*p^m)).bUnion
  (λ a : zmod (d * p ^ m), set.to_finset (equi_class m (m + 1)) a) 
\end{lstlisting}
This lemma says that any element of $\mathbb{Z}/dp^{m + 1} \mathbb{Z}$ comes from \lean{equi\_class m (m + 1) b} for some $b \in \mathbb{Z}/dp^m \mathbb{Z}$. 
The proof is now complete with the following lemma :
\begin{lstlisting}
lemma bernoulli_distribution_sum' (x : zmod (d * p^m)) : 
  ∑ (y : zmod (d * p ^ m.succ)) in 
  (λ a : zmod (d * p ^ m), ((equi_class m.succ) a).to_finset) x,
  bernoulli_distribution p d c m.succ y = bernoulli_distribution p d c m x 
\end{lstlisting}
which says, for $x \in \mathbb{Z}/dp^m \mathbb{Z}$, $E_{c, m} (x) = \sum_{y}' E_{c, m + 1} (y)$, for \lean{y ∈ equi\_class m (m + 1) x}. 

\subsection{Integrals}
The last piece in the puzzle is the integral. We use the same notation as in the previous
section. Given a measure $\mu$, and a function $f \in LC(X, R)$, $\int f d\mu := \mu(f)$. As in Theorem 12.1 of \cite{cyc}, this can be extended to a
continuous $R$-linear map $ \int_X f d\mu : C(X, R) \xrightarrow[]{} R $. 
This follows from the fact that $LC(X, R)$ is dense in $C(X, R)$; as a result, the map from 
$LC(X, R)$ to $C(X, R)$ is \href{https://leanprover-community.github.io/mathlib_docs/topology/dense_embedding.html#dense_inducing}{\lean{dense\_inducing}}, 
that is, it has dense range and the topology on $LC(X, R)$ is 
induced from the topology on $C(X,R)$. 

The continuity of the extension of the integral follows from the fact that every measure $\mu$ is uniformly continuous : 
\begin{lstlisting}
  lemma uniform_continuous (φ : measures X A) : uniform_continuous ⇑φ 
\end{lstlisting}

\subsection{Construction}
There are several possible definitions for the $p$-adic $L$-function, the most common being a meromorphic function $L_p(s, \chi)$ on
\newline $\{ s \in \mathbb{C}_p \mid |s| < p \}$ obtained by analytic continuation, such that
$$ L_p (1 - n, \chi) = -(1 - \chi \omega^{-n}(p)p^{n - 1}) \frac{B_{n, \chi \omega^{-n}}}{n} $$
for $n \ge 1$ (Theorem 5.11, \cite{cyc}). Due to the absence of $\mathbb{C}_p$ in \lean{mathlib} at the time, 
and the difficulty of showing analytic continuity(even on paper), our definition is instead motivated by Theorem 12.2, \cite{cyc}, 
which states that, for $s \in \mathbb{Z}_p$, and Dirichlet character $\chi$ with
conductor $d p^m$, with $gcd (d, p) = 1$ and $m \ge 0$, for a choice of $c \in \mathbb{Z}$
with $gcd (c, dp) = 1$ :
\begin{equation}\label{eqn:2}
  $$ (1 - \chi(c) \langle c \rangle ^{s+1}) L_p(-s, \chi) = \int_{(\mathbb{Z}/d \mathbb{Z})^{\times} \times \mathbb{Z}_p^{\times}}
\chi \omega^{-1}(a) \langle a \rangle ^s dE_c $$
\end{equation}
where $\langle a \rangle  = \omega^{-1}(a) a$, and $b^s = exp (log_p (b))$ (the exponential and logarithm are defined in 
terms of power series expansions).

Instead of using the variable $s$ (which takes values in a subset of $\mathbb{C}_p$), we choose to use an element of the weight space, 
the set of continuous monoid homomorphisms from $\mathbb{Z}/d\mathbb{Z}^{\times} \times \mathbb{Z}_p^{\times}$ to $R$. We replace
$\langle a \rangle ^s$ with \lean{w:continuous\_monoid\_hom A}. The advantage is that our $p$-adic $L$-function can now be defined over a more general space : 
a nontrivial normed commutative complete non-Archimedean $\mathbb{Q}_p$-algebra with no zero divisors. 

Given a Dirichlet character $\chi$ of level $dp^m$ with $gcd(d, p) = 1$ and $m > 0$,
we now define the $p$-adic $L$-function to be : 
$$ L_p(w, \chi) := \int_{(\mathbb{Z}/d \mathbb{Z})^{\times} \times \mathbb{Z}_p^{\times}}
\chi \omega^{-1}(a) w dE_c $$
\begin{lstlisting}
def p_adic_L_function := 
  measure.integral (bernoulli_measure R hc hc' hd na)
  ⟨(units.coe_hom R).comp (dirichlet_char_extend p d R m hd 
  (change_level _ (χ.mul ((teichmuller_character_mod_p' p R))))) * 
  w.to_monoid_hom, cont_paLf m hd _ w⟩
\end{lstlisting}
Here, \lean{dirichlet\_char\_extend} extends $\chi$ from $(\mathbb{Z}/ dp^m \mathbb{Z})^{\times}$ to 
$(\mathbb{Z}/ d \mathbb{Z})^{\times} \times \mathbb{Z}_p^{\times}$ via the restriction map. The last term 
\lean{cont\_paLf} proves the continuity of the given function, since Lean takes an element of type \lean{C((zmod d)ˣ × ℤ\_[p]ˣ, R)}. 
We have absorbed the constant term given in the LHS of (\ref{eqn:2}). This was done because Theorem 12.2 
lets $L_p(-s, \chi)$ take values in $\mathbb{C}_p$. 
In a general ring $R$, as we have chosen, division need not exist. One would then need the factor to be a unit, which may not always happen 
(for example, consider $R = \mathbb{Q}_p$). Thus, our $p$-adic $L$-function differs from the original by a constant factor. 
This factor can be easily removed if one assumes $R$ has division. 

\section{Evaluation at negative integers}
\label{section4}
We shall now prove that our chosen definition of the $p$-adic $L$-function is equivalent to the original one, that is, 
it takes the same values at negative integers : for $n > 1$,
\begin{equation}\label{eqn:1}
  $$ L_p (1 - n, \chi) = -(1 - \chi \omega^{-n}(p)p^{n - 1}) \frac{B_{n, \chi \omega^{-n}}}{n} $$
\end{equation}
For this section, we assume that $R$ is a non-Archimedean normed commutative $\mathbb{Q}_p$-algebra, 
which is complete, nontrivial, and has no zero divisors. The scalar multiplication structure obtained from $\mathbb{Q}$ 
and $\mathbb{Q}_p$ are compatible, given by \lean{is\_scalar\_tower ℚ ℚ\_[p] R} (see Section 4.2 of \cite{DD}). 
The prime $p$ is odd, and we choose positive natural numbers $d$ and $c$ which are mutually coprime and are also coprime to 
$p$. The Dirichlet character $\chi$ has level $d p^m$, where $m$ is positive. We also assume $\chi$ is even and $d$ divides 
its conductor. Let us first explain why we need the latter condition.

\subsection{Factors of the conductor}
We explain here why we need $d$ to divide the conductor of $\chi$. In this section, we do not differentiate between the associated Dirichlet 
character and the Dirichlet character. 

Recall that $\chi \omega^{-1}$ actually denotes the Dirichlet 
character multiplication of $\chi$ and $\omega^{-1}$, as explained in Section \ref{dirchar}. 
In order to translate between sums on $\mathbb{Z}/ d p^n \mathbb{Z} ^{\times}$ and $\mathbb{Z}/ d p^n \mathbb{Z}$, one needs that, for all $x \in \mathbb{Z}/ d p^n \mathbb{Z}$ 
such that $x$ is not a unit, $\chi \omega^{-k} (x) = 0$ for all $k > 0$. This is equivalent to saying, $\forall y \in \mathbb{N}$, such that 
$gcd (y, d) \ne 1$ and $gcd (y, p) \ne 1$, $gcd (y, (\chi \omega^{-k})\texttt{.conductor}) \ne 1$. 

Given coprime natural numbers $k_1, k_2$ and a character $\psi$ of level $k_1 k_2$, one can find primitive characters 
$\psi_1$ and $\psi_2$ of levels $k_1$ and $k_2$ respectively such that $\psi = \psi_1 \psi_2$ : 
\begin{lstlisting}
lemma eq_mul_of_coprime_of_dvd_conductor {m n : ℕ} [fact (0 < m * n)] 
  (χ : dirichlet_character R (m * n)) (hχ : m | χ.conductor) 
  (hcop : m.coprime n) : ∃ (χ₁ : dirichlet_character R m) 
  (χ₂ : dirichlet_character R n), χ₁.is_primitive ∧ χ = 
  χ₁.change_level (dvd_mul_right m n) * χ₂.change_level (dvd_mul_left n m) 
\end{lstlisting}
Thus, given $k > 0$, we can find primitive characters $\chi_1$ and $\chi_2$ with conductors $z_1$ and $z_2$ such that 
$z_1 \mid d$ and $z_2 \mid p^m$ and $\chi_1 \chi_2 = \chi \omega^{-k}$. The condition that $d$ divides the conductor of $\chi$ ensures 
that $z_1 = d$. As a result, if $gcd (y, d) \ne 1$, then $gcd (y, z_1 z_2) \ne 1$, so $\chi \omega^{-k} (y) = 0$ as needed.

\subsection{Main Result}
Note that the same result holds when $\chi$ is odd or when $p = 2$, the proofs differ slightly. We shall 
skip most of the details of the proof, since these are heavily computational. We shall instead highlight the key concepts 
that are used. Our reformulation of (\ref{eqn:1}) is :
\begin{lstlisting}
theorem p_adic_L_function_eval_neg_int_new :
  (p_adic_L_function m χ c na (mul_inv_pow (n - 1))) = 
  (algebra_map ℚ R) (1 / n : ℚ) *
  (1 - (χ (zmod.unit_of_coprime c _) * 
    (mul_inv_pow n (zmod.unit_of_coprime c hc', _)))) * 
  (1 - ((asso_dirichlet_character 
    (χ.mul ((teichmuller_character_mod_p' p R)^n))) p * p^(n - 1))) * 
  (general_bernoulli_number 
    (χ.mul ((teichmuller_character_mod_p' p R)^n)) n) 
\end{lstlisting}
Here, \lean{mul\_inv\_pow} is our translation of $\langle a \rangle ^s$. \newline
The proof consists of two steps : breaking up the integral in the LHS into three sums, 
and evaluating each of these sums. This is very calculation intensive, and was the longest part of the project. 
The proof is very similar to the proof of Theorem 12.2 in \cite{cyc}. 

Since $LC((\mathbb{Z}/d \mathbb{Z})^{\times} \times \mathbb{Z}_p^{\times}, R)$ is dense in $C((\mathbb{Z}/d \mathbb{Z})^{\times} \times \mathbb{Z}_p^{\times}, R)$, 
we observe that the integral $L_p (1 - n, \chi)$ is the same as :
$$ L_p (1 - n, \chi) = \lim_{j \to \infty} \sum_{a \in (\mathbb{Z}/ d p^j \mathbb{Z})^{\times}} E_{c, j} (\chi \omega^{-1} (a) \langle a \rangle ^{n - 1}) $$
\begin{equation}\label{eqn:3}
$$ = \lim_{j \to \infty} \bigg ( \sum_{a \in (\mathbb{Z}/ d p^j \mathbb{Z})^{\times}} \chi \omega^{-n} a^{n - 1} \bigg \{ \frac{a}{d p^j} \bigg \} $$
\end{equation}
\begin{equation}\label{eqn:4}
$$\qquad - \sum_{a \in (\mathbb{Z}/ d p^j \mathbb{Z})^{\times}} \chi \omega^{-n} a^{n - 1} \bigg ( c \bigg \{ \frac{c^{-1} a}{d p^j} \bigg \} \bigg ) $$
\end{equation}
\begin{equation}\label{eqn:5}
$$\qquad + \bigg ( \frac{c - 1}{2} \bigg ) \sum_{a \in (\mathbb{Z}/ d p^j \mathbb{Z})^{\times}} \chi \omega^{-n} a^{n - 1} \bigg ) $$
\end{equation}
Going from the first equation to the second took about 600 lines of code, which can be found in \href{https://github.com/laughinggas/p-adic-L-functions/blob/main/src/neg_int_eval.lean}{\lean{neg\_int\_eval.lean}}. While the proof (on paper) is only a page long, 
this is very calculation heavy in Lean, because one needs to shift between elements coerced to different types, such as $\mathbb{Z}/ (d p^j) \mathbb{Z}$, 
$\mathbb{Z}/ d \mathbb{Z} \times \mathbb{Z}/ p^j \mathbb{Z}$, $\mathbb{Z}/ d \mathbb{Z} \times \mathbb{Z}_p$, $R$ and their units. Moreover, when each of these types occur 
as locally constant or continuous functions, one needs to separately prove that each of these functions 
is also (respectively) locally constant or continuous. Other difficulties include several different ways to obtain the same term, such as \lean{equiv.inv\_fun}, 
\lean{equiv.symm}, \lean{ring\_equiv.symm} and \lean{ring\_equiv.to\_equiv.inv\_fun}. We have constructed several lemmas to simplify traversing between these terms. 

Each of these sums are then evaluated separately. The first sum in (\ref{eqn:3}) follows from Theorem \ref{thm1}, after translations between \lean{zmod (d * p\textasciicircum n)ˣ} and \lean{finset.range (d * p\textasciicircum n)}. 
This is done by the following lemma, which says 
$$ \mathbb{Z} / d p^k \mathbb{Z} \simeq \{ x \in \mathbb{N} \mid gcd (x, d) \ne 1 \} \cup \{ x \in \mathbb{N} \mid gcd (x, p) \ne 1 \} 
\cup (\mathbb{Z} / d p^k \mathbb{Z})^{\times} $$ 
\begin{lstlisting}
lemma helper_U_3 (x : ℕ) : range (d * p^x) = 
  finite.to_finset (finite_of_finite_inter 
    (range (d * p^x)) ({x | ¬ x.coprime d})) 
  ∪ ((finite.to_finset (finite_of_finite_inter 
    (range (d * p^x)) ({x | ¬ x.coprime p}))) 
  ∪ finite.to_finset (finite_of_finite_inter 
    (range (d * p^x)) ({x | x.coprime d} ∩ {x | x.coprime p}))) 
\end{lstlisting}
Each of these are made to be a \lean{finset}, since \href{https://leanprover-community.github.io/mathlib_docs/algebra/big_operators/basic.html#finset.sum}{\lean{finset.sum}} 
requires the sum to be over a \lean{finset}.We use this lemma to break our sum over \lean{finset.range (d * p\textasciicircum n)} into units and non-units. 
The condition that $d$ divides the conductor is then used to show that the associated Dirichlet character is 0 everywhere on the non-units. 
These calculations can be found in \href{https://github.com/laughinggas/p-adic-L-functions/blob/main/src/general_bernoulli_number/lim_even_character_of_units.lean}{\lean{lim\_even\_character\_of\_units.lean}}.

Evaluating the middle sum (\ref{eqn:4}) is the most tedious. It is first broken into two sums, so that the previous result can be used. Then, a 
change of variable from $a$ to $c^{-1} a$ is applied. The variable $c$ is coerced to $\mathbb{Z}/ d p^{2k} \mathbb{Z}$, increasing 
the number of coercions significantly, thus lengthening the calculations. This can be found in \href{https://github.com/laughinggas/p-adic-L-functions/blob/main/src/sum_eval/second_sum.lean}{\lean{second\_sum.lean}}.

Finally, the last sum (\ref{eqn:5}) is 0. This is where one uses that $\chi$ is even. This follows from Theorem \ref{thm2}. 
On paper, it is a one-line proof, done by substituting $a$ in the summand with $-a$ and doing caluclations mod $p^n$. However, since we work in a more general setting, 
we must go through lengthy roundabout ways instead. 
\newline Putting these sums together concludes the proof.
\section{Conclusion}
\label{section5}
\subsection{Analysis}
We list some of the observations that arose while working on this paper. 
\newline
The tactic \lean{rw} does not always work inside sums. As a result, one must use 
the \href{https://leanprover-community.github.io/extras/conv.html}{\lean{conv}} tactic to get to the expression inside the sum. While using the \lean{conv} tactic, one is said to be working 
in \lean{conv} mode. Using the \lean{conv} tactic not only lengthens the proof, but also limits the tactics one can use; 
Another way around sums is to use \href{https://leanprover-community.github.io/mathlib_docs/tactics.html#simp_rw}{\lean{simp\_rw}}, 
however, this increases compilation time of the proof. Moreover, \lean{simp\_rw} rewrites the lemma as many times as applicable, and is an unsuitable choice if one wants to apply 
the lemma just once. 

Another recurring problem was the ratio of implicit to explicit variables. The $p$-adic $L$-function, for example, has 19 arguments, of which 7 are explicit, and $p$, 
$d$ and $R$ are implicit. Excluding $R$ often means that either Lean guesses or abstracts the correct term, 
or it asks for them explicitly. In the latter case, one also gets as additional goals all the hypotheses that are dependent on $R$ and implicit, such as \lean{normed\_comm\_ring R}. 
The other alternative is to explicitly provide terms using \lean{@}, however this leads to very large expressions. 

We also ran into some instance errors. For example, since \href{https://leanprover-community.github.io/mathlib_docs/algebra/char_zero/defs.html#char_zero}{\lean{char\_zero}} 
is a class, we would like to give the lemma \lean{char\_zero R} an \lean{instance} structure. 
However, the proof is dependent on \lean{R} having the \lean{[algebra ℚ\_[p] R]} structure. Lean would then claim that this is a dangerous instance 
(for $p$ being an explicit variable) and that $p$ is a \lean{metavariable} (for $p$ being an implicit variable). Thus, we made it a \lean{lemma} instead, and had to explicitly 
feed it into implicit arguments.

While most properties regarding Bernoulli numbers and polynomials and locally constant functions have been put into \lean{mathlib}, the rest of the work is on a private repository. The author hopes to push the 
work directly to Lean 4, once the required port is complete.

\subsection{Statistics}
Given the decentralized nature of \lean{mathlib}, it is quite difficult to calculate the number of lines of code 
already existing in \lean{mathlib} which were used in this project. When initially completed, this project had about 
15000 lines of code. A major refactor was then conducted, in an effort to reduce length of individual proofs. 
We tried to uphold the spirit of \lean{mathlib}, constructing lemmas in as much generality as possible. 
The code currently consists of 28 files and about 7500 lines, grouped into appropriate categories where possible, 
according to the sections of this paper.

\subsection{Related work} 
There are several projects that require Dirichlet characters and properties of the $p$-adic integers. 
These include the project on the formalization of Fermat's last theorem for regular primes\footnote{\url{https://github.com/leanprover-community/flt-regular}}. 
There is also an effort by Prof David Loeffler which involves formalization of the classical Dirichlet $L$-function, that is somewhat dependent on this work. 
Our work on Bernoulli numbers has been used to give a formal proof of Faulhaber's theorem.

In the future, the author hopes to be able to work on Iwasawa theory, for which the $p$-adic $L$-function is a key ingredient. She also hopes to formalize more properties of Bernoulli numbers, 
that are a fundamental component of number theory.

\bibliography{ITP-2023}

\end{document}